\documentclass{amsart}
\usepackage{amsthm, amssymb, amsmath, amscd}
\usepackage{xypic}
\title[String homology of spheres and projective spaces]{String homology of spheres and projective spaces}
\author[Craig Westerland]{Craig Westerland}
\thanks{This material is based upon work supported by the National Science Foundation under agreement No. DMS-0111298.}
\input xypic

\newtheorem{theorem}{Theorem}[section]
\newtheorem{proposition}[theorem]{Proposition}

\newtheorem{lemma}[theorem]{Lemma}

\newtheorem{corollary}[theorem]{Corollary}
\newtheorem{definition}[theorem]{Definition}

\makeatletter
\let\c@equation\c@theorem
\makeatother
\numberwithin{equation}{section}

\def\R{\ensuremath{\mathbb{R}}}
\def\F{\ensuremath{\mathbb{F}}}
\def\Z{\ensuremath{\mathbb{Z}}}

\def\B{\ensuremath{{\mathcal B}}}

\def\C{\ensuremath{\mathbb{C}}}

\def\H{\ensuremath{\mathbb{H}}}
\def\cB{\ensuremath{\check{B}}}
\def\cb{\ensuremath{\check{b}}}

\def\Hom{{\rm Hom}}

\def\Tor{{\rm Tor}}
\def\Tot{{\rm Tot \,}}

\begin{document}
\bibliographystyle{amsalpha}
\begin{abstract}
We study a spectral sequence that computes the (mod $2$) $S^1$-equivariant homology of the free loop space $LM$ of a manifold $M$ (the \emph{string homology} of $M$).  Using it and knowledge of the string topology operations on $H_*(LM)$, we compute the string homology of $M$ when $M$ is a sphere or a projective space.
\end{abstract}
\maketitle
%
%
%
%
%
\section{Introduction}

The free loop space $LM = Map(S^1, M)$ of a closed $n$-manifold, $M$, admits an action of the circle $S^1$ by rotation of loops.  This space and the homotopy orbit space (or Borel construction)
$$LM_{hS^1} = LM \times_{S^1} ES^1$$
were shown in Chas and Sullivan's article \cite{cs} on string topology to admit remarkable multiplicative structures inspired by conformal field theory.  Furthermore, both spaces are entwined in the definition of topological cyclic homology given by B\"okstedt, Hsiang, and Madsen \cite{bhm}. 

The goal of this paper is to compute the homology of these spaces for certain manifolds, namely spheres and projective spaces.  Machinery for computing these (co-)homologies for general spaces does exist in the literature (see, for instance, \cite{hess, bostring, chen}).  Our purpose is to explore the link with string topology and to illustrate the power of the string topology operations in making these computations ``barehands.''  

The central idea is that the homology of $LM_{hS^1}$ may be computed via a spectral sequence (essentially Connes' spectral sequence for cyclic homology) if one has knowledge of how the homology of $S^1$ (i.e., the Batalin-Vilkovisky operator $\Delta$) acts on $H_*(LM)$.  In turn, one may often compute $\Delta$ if one understands the other (Gerstenhaber) string topology operations.

In \cite{dlops} we have computed $H_*(LM)$ when $M$ is a sphere or real, complex, or quaternionic projective space.  Furthermore, we computed certain homology operations: the one relevant to this paper is the Browder operation, a Lie bracket that we will denote $[ \cdot, \cdot ]$.  This arises in the presence of an action of the little disks operad $C_2$ on a space (or chain complex).  In the case at hand, this action is given by McClure and Smith's proof \cite{ms} of Deligne's conjecture.  A cyclic version of Deligne's conjecture \cite{kauf, msintro, tz} allows us to relate the Browder operation with the $S^1$ action, through the BV formula
$$\Delta(xy) = \Delta(x) y + (-1)^{|x|} x \Delta(y) + (-1)^{|x|} [x, y]$$
Therefore if we know the value of $\Delta$ and $[\cdot, \cdot]$ on the generators of the algebra, we may compute $\Delta$ for every class in $H_*(LM)$.  As a result, we obtain the $E_2$ term of Connes' spectral sequence for these manifolds.  A simple argument shows that for the manifolds under consideration the spectral sequence collapses at $E_2$.

In what follows, $K$ denotes one of the division algebras $\R$, $\C$, or $\H$, and $d = \dim_{\R}(K)$.  For brevity, define
$$\alpha_{d, n}(t) := t^{-d-1} + \frac{t^{d(2n)-3}}{1-t^2} \;\;\; , \;\;\; \beta_{d, n}(t):= t^{d(2n+1)-3} + \frac{1 + t^{d(2n)-2} + t^{-1}}{1-t^2}$$

\begin{theorem} \label{mainthm}

A computation of the Poincar\'e series of $H^{S^1}_*(LM; \F_2)$:

\begin{enumerate}

\item If $M=S^k$ and $k>1$, the Poincar\'e series is
$$\left(\frac{1}{1-t^{2(k-1)}} \right) \left(t^{k-1} + \frac{1+ t^{2k-1}}{1-t^2} \right)$$

\item If $M = KP^{2n+1}$ and $n>0$ if $K=\R$, the Poincar\'e series is
$$\left( \frac{t^{d(2n+1)}}{1-t^{d(2n+2)-2}} \right) \left( \frac{1-t^{-2d(n+1)}}{1-t^{-2d}} \right) \left( t^{-1} + \frac{t^{d-1} + t^{-d}}{1-t^2} \right) $$

\item If $M = KP^{2n}$, the Poincar\'e series is
$$\left( \frac{t^{d(2n)}}{1-t^{2d(2n+1)-4}} \right) \left( \left( \frac{1-t^{-2dn}}{1-t^{-2d}} \right) \alpha_{d, n}(t) + \left( \frac{1-t^{-2d(n+1)}}{1-t^{-2d}} \right) \beta_{d, n}(t) \right)$$

\end{enumerate}

\end{theorem}

We note that the Poincar\'e series for spheres agrees with answers obtained either through Carlsson-Cohen's splitting in \cite{carlcoh} or the spectral sequence defined by B\"okstedt and Ottosen in \cite{bostring}.


In section \ref{hcfsection} we introduce a cohomology theory for Frobenius algebras that is essentially the dual of negative cyclic homology.  This allows us to dualize Jones' theorem \cite{jones} identifying $H^*(LX_{hS^1})$ with the negative cyclic homology of $C^*(X)$.  A defect of this result is its limited application to formal manifolds.  Still it allows a connection with Deligne's conjecture and string topology which is employed throughout this article.

In section \ref{sssection} we use this construction to give a spectral sequence converging to $H_*(LM_{hS^1})$.  This is essentially the Bousfield-Kan spectral sequence, but its origin in cyclic cohomology allows tighter control on the differentials.  This allows us in section \ref{collapsesection} to show that the spectral sequence collapses at the $E_2$ term.  In section \ref{e2section}, we compute that $E_2$ term, proving Theorem \ref{mainthm}.

I would like to thank Igor Kriz and Ralph Kaufmann for several helpful conversations on this material.

\section{Cyclic Frobenius cohomology} \label{hcfsection}

In this section we will introduce a cohomology theory for associative Frobenius algebras $A$ which we call {\em cyclic Frobenius cohomology}, $HC_F^*(A)$.  It bears the same relation to cyclic homology that Hochschild cohomology does to Hochschild homology.  We follow Kaufmann's definition \cite{kauf} of a version (which we call $\cB$) of the $B$ operator for Hochschild \emph{cohomology} which is incorporated in the definition of $HC_F^*(A)$.  The hypothesis that $A$ is a Frobenius algebra is required in order to define $\cB$ by dualizing the usual definition of $B$ in Hochschild homology.  While we will see that $HC_F^*(A)$ encodes little more information than $HC^-_*(A)$, its virtue for our purposes will be in computing the homology of $LM_{hS^1}$.

For our purposes, a Frobenius algebra is an associative, unital, finite dimensional graded algebra $A$ over a ring $k$, endowed with a non-degenerate inner product $< \cdot, \cdot >$ which is symmetric and invariant:
$$
\begin{array}{ccc}
<a, b> = (-1)^{|a||b|} <b, a>  & , &  <ab, c> = <a, bc>
\end{array}
$$

Recall the Hochschild chain and cochain complexes of $A$:
$$CH_n(A, A) = A^{\otimes n+1} \, , \, CH^n(A, A) = Hom_k (A^{\otimes n}, A)$$
with differentials $b, \cb$ defined in the usual fashion (see, e.g., \cite{loday}).  The inner product on $A$ specifies an isomorphism between $A$ and its dual $A^*$.  We employ this as follows:  For $f \in CH^n(A, A)$, define $\tilde{f} \in \Hom_k(A^{\otimes n+1}, k) = CH_n(A, A)^*$ as
$$\tilde{f}(a_0 \otimes \cdots \otimes a_n) = <a_0, f(a_1 \otimes \cdots \otimes a_n)>$$
Since $<\cdot, \cdot>$ is nondegenerate, $f \mapsto \tilde{f}$ is an isomorphism.

Using the operator $B$ from cyclic homology, we may define an adjoint operator $\cB: CH^n(A, A) \rightarrow CH^{n-1}(A, A)$ by
$$\widetilde{\cB(f)} = B^*(\tilde{f})$$
Here $B^*$ is the linear dual of $B$.  More explicitly, $\widetilde{\cB(f)}(a_0 \otimes \cdots \otimes a_{n-1})$ is
$$\sum_{i=0}^{n-1} (-1)^{(n-1)i}<1, f(a_{n-i} \otimes \cdots \otimes a_{n-1} \otimes a_0 \otimes \cdots \otimes a_{n-i-1})> +$$
$$\sum_{i=0}^{n-1} (-1)^{(n-1)(i+1)}<a_{n-i-1}, f(1 \otimes a_{n-i} \otimes \cdots \otimes a_{n-1} \otimes a_0 \otimes \cdots \otimes a_{n-i-2})>
$$

Following the ``$B, \, b$'' definition of cyclic homology, we define a (homological) bicomplex $\B^{*, *}(A)$, the homology of whose total complex will be $HC_F^*(A)$:

\begin{definition}

Define the bicomplex $\B^{*, *}(A)$ using the Hochschild cochain complex $CH^*(A, A)$:
$$\B^{p, q}(A) := CH^{p-q}(A, A) = \Hom_k(A^{\otimes p-q}, A)$$
for $p-q \geq 0$ and $p \geq 0$.  The vertical differential $\cb: \B^{p, q}(A) \rightarrow \B^{p, q-1}(A)$ is the Hochschild cohomology differential.  The horizontal differential $\cB: \B^{p, q}(A) \rightarrow \B^{p-1, q}(A)$ is defined above.

\end{definition}

By comparison with the same result from cyclic homology, one can show the following (this is where the assumptions of symmetry and invariance are used):

\begin{lemma}

$\B^{*, *}(A)$ is a bicomplex; that is, $\cB^2 = 0 = \cb^2$ and $\cB \cb+ \cb \cB = 0$.

\end{lemma}

\begin{definition}

The {\em cyclic Frobenius cochain complex of $A$} is defined to be
$$CC_F^*(A) := \Tot (\B^{*, *}(A))$$
and its homology, $HC_F^*(A)$, is the {\em cyclic Frobenius cohomology of $A$}.

\end{definition}

Poincar\'e duality is used to prove the following classical fact about our main example.

\begin{proposition} \label{poincareprop}

The cohomology algebra $H^*(M; k)$ of a $k$-oriented closed $n$-manifold $M$ is a Frobenius algebra.

\end{proposition}

\begin{proof}

Here we grade $H^*(M)$ negatively to give it a homological differential.  $H^*(M)$ admits a graded commutative, associative cup product $\smile$.  The inner product is the intersection form: Let $[M] \in H_n(M)$ be the fundamental class of $M$.  Then for $a, b \in H^*(M)$, the inner product is defined to be the evaluation of the cup product on the fundamental class:
$$<a, b> := (a \smile b)([M])$$
That $< \cdot, \cdot>$ is graded symmetric and invariant follows from the graded commutativity and associativity of the cup product.

To show that the inner product is nondegenerate, for each $a \in H^k(M)$, we must produce a class $b \in H^{n-k}(M)$ for which $<b, a> \neq 0$.  Choose $b$ to be any cohomology class which is nonzero on the homology class
$$a \frown [M] \in H_{n-k}(M)$$
Then 
$$<b, a> = (b \smile a)([M]) = b(a \frown [M]) \neq 0$$

\end{proof}

\begin{theorem} \label{cyclicstringtheorem}

If $M$ is formal, the cyclic Frobenius cohomology of its cohomology algebra is isomorphic to its string homology:

$$HC_F^*(H^*(M)) \cong \Sigma^{-n} H^{S^1}_*(LM)$$

\end{theorem}

One would like an version of this result without an appeal to formality.  It seems clear that for such a result, one needs to replace $H^*(M)$ with some version of the cochain complex of $M$.  However, it is not apparent to us how to endow $C^*(M)$ with the structure of a Frobenius algebra.  One can, for instance, triangulate $M$ and use the simplicial cochain algebra, with the same inner product as described in Proposition \ref{poincareprop}.  Unfortunately, on the whole of the cochain complex, this inner product is degenerate.  This suggests the need for a notion of a {\em homotopy Frobenius algebra} in which a version of $C^*(M)$ would be a prime example, and for which it would be possible to define cyclic Frobenius cohomology.  For our purposes we shall only be considering formal manifolds, and therefore will not explore such subtleties.

We refer the reader to the work of Xiaojun Chen \cite{chen} where a very similar model for the $S^1$-equivariant chain complex of $LM$ that may avoid such difficulties is developed using methods of rational homotopy theory and Brown's twisting cochains.

\begin{proof}

Jones has shown in \cite{jones} that 
$$HC^-_*(C^*(X)) \cong H_{S^1}^*(LX)$$
for any space $X$.  Here $HC^-_*$ is the negative cyclic homology functor.  If $X$ is formal, we may obviously replace $C^*(X)$ with $H^*(X)$ in the isomorphism above.  Therefore there is a quasi-isomorphism
$$CC_*^-(H^*(X)) \simeq C^*(LX \times_{S^1} ES^1)$$
where $CC^-_*$ is the chain complex which computes cyclic homology.  $CC^-_*(A)$ is the totalization of a bicomplex $\B^-_{*, *}(A)$, with
$$\B^-_{p, q} (A):= CH_{q-p} (A, A) \;\; , \;\; p \leq 0$$
The vertical differential is the Hochschild homology differential $b$, and the horizontal differential is the operator $B$.

The theorem will follow if we exhibit a (degree shifting) isomorphism of bicomplexes
$$c: \B^{*, *} (H^*(M)) \rightarrow (\B^-_{*, *}(H^*(M)))^*$$
for manifolds $M$, since the totalization of the latter bicomplex is equivalent to $C_*(LM \times_{S^1} ES^1)$.  This isomorphism will negate degrees: for $p \geq 0$
$$c: \B^{p, q} (H^*(M)) \rightarrow (\B^-_{-p, -q}(H^*(M)))^*$$
The domain is $CH^{p-q}(H^*(M), H^*(M))$ and the range is $CH_{p-q}(H^*(M), H^*(M))^*$.  So we may define $c(f):= \tilde{f}$.

It is definitional that $c \circ \cB = B^* \circ c$, and a computation (again relying on symmetry and invariance) that $c \circ \cb = b^* \circ c$.  So $c$ is a map of bicomplexes.  We have already seen that it is an isomorphism by the nondegeneracy of the inner product.  The theorem follows.

\end{proof}

It is worth pointing out that the above proof, once divorced from topological applications, gives an isomorphism $\B^{*, *} (A) \cong (\B^-_{*, *}(A))^*$ for any Frobenius algebra $A$.

\section{Connes' spectral sequence} \label{sssection}

%
%
%
%
%
%
%
%
%
%
%
%
%
%
%

In this section, we introduce a spectral sequence coverging to $HC_F^*(A)$ and when $A = H^*(M)$, we relate it to the Bousfield-Kan spectral sequence for the simplicial space $LM_{hS^1}$.  Write the fundamental class of $S^1$ as $\Delta$ so that as a ring the homology of $S^1$ is an exterior algebra on $\Delta$:
$$H_*(S^1) = \Lambda[\Delta]$$
and $H_*(LM)$ becomes a module over this ring through the action of $S^1$ on $LM$.

One may filter $CC^*(A)$ by vertical stripes in the bicomplex $\B^{*, *}(A)$.  This, in turn, produces a spectral sequence converging to $HC_F^*(A)$; the $E_1$-term of this spectral sequence is given by
$$E_1^{p, q} = HH^{p-q}(A, A); \; p \geq 0$$
with differential $d_1$ (of bidegree $(-1, 0)$) given by $B$.  In the dual case, the analogous spectral sequence for cyclic \emph{homology} was considered by Connes and called \emph{Connes' spectral sequence} in \cite{weibel}.  We keep that terminology here.


When $M$ is simply connected (and in certain other cases such as $\R P^n$; see \cite{dlops}), the homology of $LM$ can be computed using the results of \cite{cj}.  In that article, Cohen and Jones introduce a spectrum $LM^{-TM}$ (the Thom spectrum of the pullback of $-TM$ via the evaluation at $1$, $ev:LM \rightarrow M$), and show that
$$H_*(LM^{-TM}) \cong HH^*(C^*(M), C^*(M))$$
Via the Thom isomorphism, we may identify $H_*(LM)$ as equivalent to the $n$-fold suspension of $HH^*(C^*(M), C^*(M))$. 

So, when $M$ is formal, the $p^{\rm th}$ column in the $E_1$ term of Connes' spectral sequence for $HC_F^*(H^*(M))$ is isomorphic to the $p^{\rm th}$ suspension of $\Sigma^{-n}H_*(LM)$.  The differential $d_1$, induced by $\cB$, is a map between the columns.  The following is an application of Theorem 4.1 of \cite{jones}, dualized as in the arguments presented in the proof of Theorem \ref{cyclicstringtheorem}.

\begin{lemma}

For $x$ in the $E_1$ term of Connes' spectral sequence, $d_1(x) = \cB(x) = \Delta(x)$.

\end{lemma}

This implies that Connes' spectral sequence is essentially the Bousfield-Kan spectral sequence for $X_{hS^1}$ after $E_2$.  This is because the $E_2$-term of that spectral sequence is $\Tor_{H_*(S^1)}^*(\Z, H_*(LM))$.  Computing this using the standard periodic resolution of $\Z$ over $H_*(S^1) = \Lambda [ \Delta ]$ gives the $E_1$ term of Connes' spectral sequence.

We collect information about the differentials in the spectral sequence that will allow us to prove that it collapses for the manifolds under consideration.  The first statement below is standard; the second follows from the fact that $\cB$ raises topological degree by one.

\begin{lemma} \label{differentialcorollary}

The $r^{\rm th}$ differential in Connes' spectral sequence for $HC_F^*(H^*(M))$ is of bidegree $(-r, r-1)$.  That is,
$$d_r: E_r^{p, q} \rightarrow E_r^{p-r, q+r-1}$$
is a map between a subquotient of $HH^{p-q}(H^*(M), H^*(M))$ and a subquotient of $HH^{p-q+1-2r}(H^*(M), H^*(M))$.  Moreover, $d_r$ is of topological degree $+1$ as a map between subquotients of $\Sigma^{-n} H_*(LM)$.

\end{lemma}

\section{A computation of the $E_2$ term of Connes' spectral sequence} \label{e2section}

To compute the $E_2$ term of Connes' spectral sequence, we determine the action of the operator $\Delta$ on $HH^*(C^*(M), C^*(M))$.  This is accomplished for generators through direct computations or filtration arguments.  It is extended by the BV formula
$$\Delta(xy) = \Delta(x) y + (-1)^{|x|} x \Delta(y) + (-1)^{|x|} [x, y]$$
and a computation of the bracket.  We recall from \cite{dlops} the following computations:

\begin{enumerate}

\item If $k>1$, $HH^*(C^*(S^k), C^*(S^k))$ is isomorphic as an algebra to $\F_2[x, v]/(x^2)$, where the dimensions of $x$ and $v$ are $-k$ and $k-1$, respectively. The bracket is given by $[x, v] = 1$.

\item Let $K$ be one of $\R$, $\C$, or $\H$, and let $d = \dim_{\R}(K)$. For $n$ odd (and greater than $1$ if $K=\R$),
$$HH^*(C^*(K P^n), C^*(K P^n)) = \F_2[x, v, t]/(x^{n+1}, v^2 - \frac{n+1}{2} t x^{n-1})$$
and for $n$ even,
$$HH^*(C^*(K P^n), C^*(K P^n)) = \F_2[x, u, t]/(x^{n+1}, u^2, t x^n, u x^n)$$
where the topological dimensions of $x, u, v$, and $t$ are $-d$, $-1$, $d-1$, and $d(n+1)-2$ respectively.  Their Hochschild degrees are $0, 1, 1$, and $2$.  The bracket is given on generators by
$$
\begin{array}{lllll}
[x, v] = 1, & [x, u] = x, & [x, t] = 0, & [v, t] =0, & [u, t] = t.
\end{array}
$$

\end{enumerate}

In general, we will use the notation $M$ to refer to any of the manifolds $S^k$ ($k>1$), $\R P^n$ ($n>1$), $\C P^n$, or $\H P^n$.

\begin{lemma} \label{generatorlemma}

$\Delta$ vanishes on algebra generators of $HH^*(C^*(M), C^*(M))$.

\end{lemma}

\begin{corollary} \label{deltalemma}

\noindent \begin{enumerate}

\item In $HH^*(C^*(S^k), C^*(S^k))$, $\Delta(x^a v^b) = ab x^{a-1} v^{b-1}$.

\item  Depending upon the parity of $n$, monomials in $HH^*(C^*(KP^n), C^*(KP^n))$ may be written as $x^a v^b t^c$ or $x^a u^b t^c$ (where $b = 0, 1$).  Then
$$
\begin{array}{lcr}
\Delta(x^a v^b t^c) = ab x^{a-1} t^c & and & \Delta(x^a u^b t^c) =  (a+c)b x^a t^c
\end{array}
$$

\end{enumerate}

\end{corollary}

\begin{proof}

Assuming Lemma \ref{generatorlemma} we will prove part (1).  Part (2) is somewhat tedious and proved in the same fashion.

First notice that from the Leibniz formula
$$[\alpha, \beta \gamma] = [\alpha, \beta] \gamma + [\alpha, \gamma] \beta,$$
it follows that
$$[\alpha, \beta^p] = p[\alpha, \beta] \beta^{p-1}. \leqno{(*)}$$

So taking $a=0$ in part (1), we know that
$$\Delta(v^b) =  \Delta(v) v^{b-1} + v \Delta(v^{b-1}) + [v, v^{b-1}]$$
The first term is $0$ by Lemma \ref{generatorlemma}, and the third by $(*)$.  Therefore $\Delta(v^b) = 0$ by induction.  So
$$
\begin{array}{rcl}
\Delta(xv^b) & = & \Delta(x) v^b + x \Delta(v^b) + [x, v^b] \\
             & = & b [x, v] v^{b-1} \\
             & = & bv^{b-1}
\end{array}
$$
Part (1) follows.

\end{proof}

\noindent {\it Proof of Lemma \ref{generatorlemma}.}

The operator $\Delta$ (induced by $\cB$) lowers the Hochschild degree by $1$: if $\alpha \in HH^p(R, R)$, $\Delta (\alpha) \in HH^{p-1}(R, R)$.  Automatically, we thereby obtain
$$\Delta(x) = 0,$$
since $x \in HH^0(C^*(M), C^*(M))$.

The element $v$ has Hochschild degree $1$ and topological degree $d-1$, where $d$ is as above if $M=KP^n$ and $d = k$ if $M = S^k$.  Therefore $\Delta(v)$ has Hochschild degree $0$ and topological degree $d>0$.  Since there are no elements in $HH^0(C^*(M), C^*(M))$ of positive topological degree, $\Delta(v) = 0$.

Similarly, $\Delta(t) \in HH^1(C^*(KP^n), C^*(KP^n))$ has topological degree $d(n+1)-1$.  If nonzero, $\Delta(t)$ may be written as
$$
\Delta(t) = \left\{
\begin{array}{cc}
x^k v; & n \; {\rm is \; odd} \\
x^k u; & n \; {\rm is \; even}
\end{array} \right.
$$
In the first case, the topological degree of $x^k v$ is $-kd +d-1$, so we must have $k = -n <0$, which is impossible.  Similarly, if $n$ is even, $k = -n-1 <0$.

Finally, to show that $\Delta(u) = 0$ we use the description of $\Delta$ as induced by the $\cB$ operator.  In \cite{dlops} we found that a representative for the class $u$ is the function $\overline{u} \in CH^1(\F_2[x]/x^{n+1}, \F_2[x]/x^{n+1})$ given by
$$\overline{u}: x^m \mapsto m x^m$$
so
$$
\begin{array}{rcl}
\widetilde{\cB(\overline{u})}(x^m) & = & <1, \overline{u}(x^m)> + <x^m, \overline{u}(1)> \\
 & = & m < 1, x^m> + <x^m, 0>\\
 & = & m \cdot x^m([KP^n])
\end{array}$$
For $x^m([KP^n])$ to be nonzero, $m = n$.  Since $n$ is even, the product is zero.

\qed

\noindent {\it Proof of Theorem \ref{mainthm}.}

Using Corollary \ref{deltalemma} we will compute the Poincar\'e series of the $E_2$ term of Connes' spectral sequence.  Lemma \ref{collapselemma} then gives us Theorem \ref{mainthm}.  For brevity, we only do this computation for $M = KP^{2n+1}$.  The computations for other manifolds are similar; the case for spheres is easier, the case for even projective spaces is more tedious.

Recall that $E_1^{p, q} = HH^{p-q}(H^*(M), H^*(M))$, and that through Cohen and Jones' work, the $p^{\rm th}$ column of the $E_1$ term of the spectral sequence may be thought of as $\Sigma^{p-n} H_*(LM)$.  Therefore we will adopt a new notation for the spectral sequence where the second variable indicates the topological degree, rather than the Hochschild degree:
$$E_1^{p, q} = H_{q+n-p}(LM)$$
More carefully, we are employing the fact that our spectral sequence is actually triply graded -- by filtration degree, Hochschild degree, and topological degree -- and within a fixed filtration grading (the columns), we re-grade using the topological degree.  In this format, the spectral sequence is graded identically to the Bousfield-Kan spectral sequence for $LM_{hS^1}$.

Examine the action of the operator $\Delta$ on $H_*(LKP^{2n+1})$.  There are three types of classes:

\begin{enumerate}

\item Classes $a$ for which $\Delta(a) \neq 0$.  We say these ``survive alone.''
\item Classes $b$ for which $\Delta(b) = 0$ and there is a class $b'$ with $\Delta(b') = b$.  We say these are ``hit.''

\item Classes $c$ for which $\Delta(c) = 0$ that are not in the image of $\Delta$.  We say these classes ``propagate a stripe.''

\end{enumerate}

Recall that the $E_1$ term of the spectral sequence looks like
$$
\xymatrix@1{ \Sigma^{-n} H_*(LM) & \Sigma^{1-n} H_*(LM) \ar[l]_-{\Delta} & \Sigma^{2-n} H_*(LM) \ar[l]_-{\Delta} & \cdots \ar[l]_-{\Delta} }
$$
(where each $H_*(LM)$ is a column).  

Classes $a$ which survive alone give rise to an element of $E_2^{0, |a|}$.  Classes $b$ which are hit do not give rise to any element of $E_2$.  Classes $c$ which propagate a stripe give a class in $E_2^{p, |c|+p}$ for each $p \geq 0$.

Examining Corollary \ref{deltalemma}, we see that for every $k$ and $c$, $x^{2k+1}vt^c$ survives alone, $x^{2k}vt^c$, and $x^{2k+1}t^c$ propagate a stripe, and all other monomials are hit.

The dimension of $x^{2k+1}vt^c$ is 
$$|x^{2k+1}vt^c| = -1 + k(-2d) + c(d(2n+2)-2)$$
so the Poincar\'e series of the space that they span is
$$t^{-1} \sum_{k=0}^n \sum_{c=0}^{\infty} t^{k(-2d)} t^{c(d(2n+2) -2)} = t^{-1} \left( \frac{1-t^{-2d(n+1)}}{1-t^{-2d}} \right) \left( \frac{1}{1 - t^{d(2n+2)-2}} \right)$$
This is exactly the contribution to $E_2^{0, *}$ of the elements that survive alone.

The dimension of $x^{2k}vt^c$ is
$$|x^{2k}vt^c| = (2k)(-d) + d-1 + c(d(2n+2)-2)$$
Similarly, the Poincar\'e series of the family $<x^{2k}vt^c>$ is
$$t^{d-1} \left( \frac{1-t^{-2d(n+1)}}{1-t^{-2d}} \right) \left( \frac{1}{1 - t^{d(2n+2)-2}} \right),$$
and for the family $<x^{2k+1}t^c>$:
$$t^{-d} \left( \frac{1-t^{-2d(n+1)}}{1-t^{-2d}} \right) \left( \frac{1}{1 - t^{d(2n+2)-2}} \right)$$

Each element of dimension $q$ in the latter two families gives rise to a sequence of classes in $E_2^{p, q+p}$ for each $p$.  So in counting their contribution to $E_2$, we must multiply the answer by
$$\sum_{k=0}^{\infty} t^{2k} = \frac{1}{1-t^2}$$

Adding these three series gives us the Poincar\'e series for $E_2$:
$$\left(t^{-1} + \frac{t^{d-1} + t^{-d}}{1-t^2} \right) \left( \frac{1-t^{-2d(n+1)}}{1-t^{-2d}} \right) \left( \frac{1}{1 - t^{d(2n+2)-2}} \right)$$

Recall that the spectral sequence computes a {\em desuspension} of the homology of $LM$; one needs to multiply this series by $t^{\dim M} = t^{d(2n+1)}$ to get the correct answer.

\qed

%
%
%
%
%
%
%
%
%
%
%
%
%
%
%
%
%
%

\section{Collapse of the spectral sequence} \label{collapsesection}

We complete the proof of Theorem \ref{mainthm} with the following result:

\begin{lemma} \label{collapselemma}

For the manifolds $M$ considered in this paper, Connes' spectral sequence for $HC_F^*(C^*(M))$ collapses at the $E_2$-term.

\end{lemma}

\begin{proof}

From Corollary \ref{differentialcorollary} we know that all differentials $d_r$ in the spectral sequence are of the form
$$
\begin{array}{ccc}
{\rm subquotient} & & {\rm subquotient} \\
{\rm of} & \longrightarrow & {\rm of} \\
HH^{k}(C^*(M), C^*(M)) & & HH^{k+1-2r}(C^*(M), C^*(M))
\end{array}
$$
of topological degree $+1$.  We will show that if $r > 1$, such a map is $0$ by examining the range of topological degrees of the source and target.  For simplicity, we take $M=KP^n$ with $n$ odd; the proofs for even projective spaces and spheres are similar.

If $k = 2l$ is even, then we showed in \cite{dlops} that $HH^{k}(C^*(M), C^*(M))$ is concentrated in topological degrees 
$$l(d(n+1)-2) - jd; \; j \in \{ 0, \ldots, n \}$$
and if $k = 2l+1$ is odd, it lies in dimensions
$$l(d(n+1)-2) + d-1 - jd; \; j \in \{ 0, \ldots, n \}$$

Consider $d_r$ as a mapping from a subquotient of $HH^{2l}$ to a subquotient of $HH^{2l+1-2r}$, and let $\alpha$ lie in the domain.  The smallest possible dimension for $\alpha$ is 
$$l(d(n+1)-2) - nd = (l-1)nd + ld -2l$$ 
(realized by the element $x^n t^l$) and therefore, the dimension of $d_r(\alpha)$ is
$$|d_r(\alpha)| = |\alpha|+1 \geq (l-1)nd + ld -2l +1$$
We claim that this is larger than the dimension of any element in the range.  The largest possible dimension in the range is
$$\begin{array}{lll}
(l-r)(d(n+1)-2) + d-1 & = & lnd + ld -2l + (-rnd -rd + 2r + d-1) \\
                      & = & (l-1)nd + ld -2l +1 + (r-1)(2-d-nd)
\end{array}
$$
(realized by the element $v t^{l-r}$).  Since $d, n \geq 1$, but not $d = n = 1$ (in which case we would be considering $\R P^1 = S^1$), and $r>1$,
$$(r-1)(2-d-nd) < 0$$
so $|d_r(\alpha)| > |\beta|$ for every $\alpha \in HH^{2l}$ and $\beta \in HH^{2l+1-2r}$.

To check that $d_r$ is $0$ as a map from a subquotient of $HH^{2l+1}$ to a subquotient of $HH^{2l+2-2r}$ takes only a little more work.  First we notice that the element of the domain of lowest dimension, $x^n v t^l$ does not lie in $E_2$, since
$$\Delta (x^n v t^l) = x^{n-1} t^l \neq 0$$
so the smallest dimension of an element $\alpha$ for which $d_r(\alpha)$ might be nonzero is
$$l(d(n+1)-2) + d-1 - d(n-1)$$
corresponding to $x^{n-1} v t^l$.  The class of largest possible dimension in the range of $d_r$ is $t^{l+1-r}$, of dimension $(l+1-r)(d(n+1)-2)$.  The difference in dimension between $d_r(\alpha)$ and the largest possible target is then
$$|t^{l+1-r}| - |d_r(\alpha)| \leq |t^{l+1-r}| - |d_r(x^{n-1} v t^l)| = (2dn-d-2) - r(dn+d-2).$$
Since $r\geq 2$,
$$\begin{array}{lll}
|t^{l+1-r}| - |d_r(\alpha)| &   =  & (2dn-d-2) - r(dn+d-2) \\
                            & \leq & (2dn-d-2) - 2(dn+d-2) \\
                            &   =  & 2-3d \\
                            &   <  & 0
\end{array}
$$
since $d>1$.  So again, $|d_r(\alpha)| > |\beta|$ for every $\alpha \in HH^{2l+1}$ and $\beta \in HH^{2l+2-2r}$.

\end{proof}

\bibliography{biblio}

\end{document}